\documentclass[11pt,twoside,leqno]{aomamlt2e} 
  \pageno{1377}
\received{May 24, 2004}

 \renewcommand{\thesection}{\arabic{section}}
\makeatletter
\renewcommand{\@seccntformat}[1]{\csname
the#1\endcsname.\hspace{0.5em}\setcounter{Subsec}{0}\setcounter{Subsubsec}{0}}\makeatother

\newtheorem{thm}{Theorem}[section]

\newtheorem{lemma}[thm]{Lemma}
\theorembodyfont{\upshape}

\newtheorem{rem}[thm]{{\it Remark}}

\newcommand{\bb}{\begin{equation}}
\newcommand{\ee}{\end{equation}}
\newcommand{\bq}{\begin{eqnarray}}
\newcommand{\eq}{\end{eqnarray}}
\newcommand{\bqn}{\begin{eqnarray*}}
\newcommand{\eqn}{\end{eqnarray*}}

\def\theequation{\thesection.\arabic{equation}}

\begin{document}
\currannalsline{162}{2005} 

 \title{Formation of singularities for a transport\\ equation with nonlocal velocity}

 \acknowledgements{Partially supported by BFM2002-02269 grant.\hfill\break
\hglue12pt$^{\textstyle\ast\ast}$Partially supported by BFM2002-02042 grant.\hfill\break
\hglue3pt$^{\textstyle\ast\ast\ast}$Partially supported by BFM2002-02042 grant.}
 
\twoauthors{Antonio C\'ordoba$^{\textstyle\ast}$,  Diego C\'ordoba$^{\textstyle\ast\ast}$,}{Marco A.
Fontelos\lower.75pt\hbox{$^{\textstyle\ast\ast}$}}

 \institution{Universidad Aut\'{o}noma de Madrid, Madrid, Spain\\
\email{antonio.cordoba@uam.es}\\
\vglue-9pt
Consejo Superior de Investigaciones Cient\'{i}ficas, 
   Madrid, Spain\\
\email{dcg@imaff.cfmac.csic.es}\\
\vglue-9pt
 Universidad Rey Juan Carlos, Madrid, Spain
\\
\email{mafontel@escet.urjc.es}}

 \shorttitle{Formation of singularities for a transport equation}
  
  \shortname{Antonio C\'ordoba,  Diego C\'ordoba, and Marco A.
Fontelos}

\centerline{\bf Abstract}
\vskip12pt
We study a 1D transport equation with nonlocal velocity and show
 the formation of singularities in finite time for a generic family of initial data.
  By adding a diffusion term the finite time singularity is prevented and
   the solutions exist globally in time.

\section{Introduction}

  In this paper we study the nature of the solutions to the
  following class of equations
  \begin{eqnarray}
\theta _{t}-\left( H\theta \right) \theta _{x}=-\nu
\Lambda^{\alpha}\theta, \quad\quad x\in R\; \label{1eq}
\end{eqnarray}
where $H\theta$ is the Hilbert transform defined by
\begin{eqnarray*}
H\theta\equiv \frac{1}{\pi} PV\int \frac{\theta(y)}{x-y}dy,
\label{hilbert}
\end{eqnarray*}
$\nu$ is a real positive number, $0\leq\alpha\leq 2$  and
$\Lambda^{\alpha}\theta\equiv(-\Delta)^{\frac{\alpha}{2}}\theta$.

 This equation represents the simplest case of a transport
equation with a nonlocal velocity and with a viscous term
involving powers of the laplacian. It is well known that the
equivalent equation with a local velocity $v=\theta$, known as
Burger«s equation, may develop shock-type singularities in finite
time when $\nu=0$ whereas the solutions remain smooth at all times
if $\nu>0$ and $\alpha=2$. Therefore a natural question to pose is
whether the solutions to (\ref{1eq}) become singular in finite
time or not depending on $\alpha$ and $\nu$. In fact this question
has been previously considered in the literature motivated by the
strong analogy with some important equations appearing in fluid
mechanics, such as the  3D Euler incompressible vorticity equation and the 
Birkhoff-Rott equation modelling the evolution of a vortex sheet,
where a crucial mathematical difficulty lies in the nonlocality of
the velocity. Since the fundamental problem concerning both 3D
Euler and Birkhoff-Rott equations is the formation of
singularities in finite time, the main goal of this paper will be
to solve this issue for the model (\ref{1eq}).

3D Euler equations, in terms of the vorticity vector are
\begin{eqnarray}
\omega _{t} + v\cdot\nabla \omega = \omega D(\omega) \label{euler}
\end{eqnarray}
where $D(\omega)$ is a singular integral operator of $\omega$
whose one dimensional analogue is the Hilbert transform and the
velocity is given by the Biot-Savart formula in terms of $\omega$.
In order to construct lower dimensional models containing some of
the main features of (\ref{euler}), Constantin, Lax and Majda
\cite{CLM} considered the scalar equation
\begin{eqnarray}
\omega _{t} + v \omega _{x}= \omega H\omega  ; \label{eqq}
\end{eqnarray}
with $v=0$ and showed existence of finite time singularities. The
effect of adding a viscous dissipation term has been studied in
\cite{Schochet}, \cite{Vasudeva}, \cite{Vaswegert}, \cite{Yang}\break
and \cite{Sakajo}. In order to incorporate  the advection term
$v\omega_x$ into the model, De Gregorio, in \cite{De1} and
\cite{De2}, proposed a velocity given by an integral operator of
$\omega$. If we take an $x$ derivative of (\ref{1eq}) and define
$\theta_x\equiv\omega$ we obtain a viscous version of the equation
(\ref{eqq}) with $v= -H\theta$ which is similar to the one
proposed in {\cite{De1} and \cite{De2}.

The analogy of (\ref{1eq}) with Birkhoff-Rott equations was first
established in \cite{B} and \cite{M}. These are
integrodifferential equations  modelling the evolution of vortex
sheets with surface tension. The system can be written in the form
\begin{eqnarray}
\frac{\partial}{\partial t}z^{*}(\alpha,t)&=& \frac{1}{2\pi i} PV
\int \frac{\tilde{\gamma}(\alpha')d\alpha'}{z(\alpha,t) - z(\alpha',t)} \label{1eq1} \\
 \frac{\partial\tilde{\gamma}}{\partial t}&=&\sigma\kappa_{\alpha} \label{1eq2}
\end{eqnarray}
where $z(\alpha,t)= x(\alpha,t) + i y(\alpha,t)$ represents the
two dimensional vortex sheet parametrized with $\alpha$, and where
$\kappa$ denotes mean curvature. Following \cite{B} we substitute,
in order to build up the model, the equation (\ref{1eq1}) by its
1D analog
\begin{eqnarray*}
\frac{dx(\alpha,t)}{dt}= -H(\theta)\label{eq3}
\end{eqnarray*}
where we have identified $\gamma(\alpha,t)$ with $\theta$. In the
limit of $\sigma =0$ in (\ref{1eq2}) we conclude that $\gamma$ is
constant along trajectories and this fact leads, in the 1D model,
to the equation
\begin{eqnarray}
\theta _{t}-\left( H\theta \right) \theta _{x}=0. \label{1eq4}
\end{eqnarray}

 There is now overwhelming evidence
that vortex sheets form curvature singularities in finite time.
This evidence comes back from the classical paper by Moore
\cite{Moore} where he studied the Fourier spectrum of
$z(\alpha,t)$ and, in particular, its asymptotic behavior when the
wave number $k$ goes to infinity. His numerical results showed
that, up to very high values of $k$, this asymptotic behavior is
compatible with the formation of a curvature singularity in finite
time. Although there has been a very intense activity in order to
provide a definitive proof of the formation of such a singularity
(see discussions and references in \cite{Moore}, \cite{BM} and
\cite{B})  the  existing results  are mostly supported in
numerics or formal asymptotics and do not constitute a full
mathematical proof. The same kind of argument was used in
\cite{B} in order to show the existence of singularities for the
1D analog (\ref{1eq4}).

The system (\ref{1eq1}) and (\ref{1eq2}) with $\sigma=0$ has the
very interesting property of being ill-posed for general initial
data. A linear analysis of small perturbations of planar sheets
leads to catastrophically growing dispersion relations. Several
attempts at regularization were introduced through the
incorporation of effects, such as surface tension or viscosity
(see \cite{BM} for a comprehensive review). In the same spirit we
will also study the effects of artificial viscosity terms on the
solutions for our model. More precisely we will prove the
existence of blow-up in finite time for (\ref{1eq}) with  $\nu=0$
in Section 2 and, conversely, the global existence of solutions
when $\nu>0$ and $1<\alpha\leq 2$ in Section 3.

\section{Blow-up for $\nu=0$}

The local existence of solutions to (\ref{1eq}) was established in
\cite{B}.  In this section we will show the existence of
blowing-up solutions to (\ref{1eq4}) for a generic class of
initial data.

Let us consider a symmetric, positive, and $C^{1+\varepsilon
}(\mathbb{R)}$ initial profile $\theta =\theta _{0}(x)$ such that
${\rm max}_x\theta_0=\theta _{0}(0)=1$. We will also assume
\begin{eqnarray*}
{\rm Supp}(\theta _{0}(x))\subset \left[ -L,L\right]\; .
\end{eqnarray*}
Under these assumptions, it is clear that $\theta (x,t)$ will
remain positive (given the transport character of equation
(\ref{1eq}) for $\nu =0$) and symmetric. Then, $H\theta $ will be
antisymmetric and positive for $x\geq L$. This implies the
following properties for $\theta(x,t)$:
\begin{eqnarray*}
 &\centerdot& {\rm Supp}(\theta (x,t))\subset \left[ -L,L\right]\;,\\
 &\centerdot& {\rm max}_x\theta=\theta (0,t)=1\; ,\\
 &\centerdot& \|\theta\|_{L^1}(t)\leq \|\theta\|_{L^1}(0)\;,\\
 &\centerdot& \|\theta\|_{L^2}(t)\leq \|\theta\|_{L^2}(0)\;.
\end{eqnarray*}

\begin{thm}
Under the conditions stated above for $\theta_0${\rm ,} the solutions of
{\rm (\ref{1eq})} with $\nu=0$ will always be such that
$\|\theta_x\|_{L^{\infty}}$ blows up in finite time.
\end{thm}

\Proof   Since $\theta _{t}=-(1-\theta )_{t}\equiv
-f_{t}$, $\theta _{x}=-(1-\theta
)_{x}\equiv -f_{x}$ and $H\theta =-H(1-\theta )\equiv -Hf$, we can write, from (\ref{1eq4}),
\begin{eqnarray}
(1-\theta )_{t}=-H(1-\theta )(1-\theta )_{x}\;.  \label{eq2}
\end{eqnarray}

 We now divide  (\ref{eq2}) by $x^{1+\delta }$ with $0<\delta <1$,
integrate in $\left[ 0,L\right] $ and obtain the following
identity:
\begin{eqnarray}
\frac{d}{dt}\left( \int_{0}^{L}\frac{(1-\theta )}{x^{1+\delta
}}dx\right) =-\int_{0}^{L}\frac{(1-\theta )_{x}H(1-\theta
)}{x^{1+\delta }}dx\; . \label{eq3}
\end{eqnarray}

Given  the fact that $\theta $ vanishes outside the interval
$\left[ -L,L\right] $, we can write the right-hand side of (\ref
{eq3}) in the form
\begin{equation}
-\int_{0}^{L}\frac{(1-\theta )_{x}H(1-\theta )}{x^{1+\delta }}%
dx=-\int_{0}^{\infty }\frac{(1-\theta )_{x}H(1-\theta
)}{x^{1+\delta }}dx\;. \label{eq4}
\end{equation}

In the next lemma we provide an estimate for the right-hand side
of (\ref {eq4}).

\begin{lemma}
Let $f\in C_{c}^{\infty }(\mathbb{R}^{+})$. Then for $0<\delta <1$
there exists a constant $C_{\delta }$ such that
\begin{eqnarray}
-\int_{0}^{\infty }\frac{f_{x}(x)(Hf)(x)}{x^{1+\delta }}dx\geq
C_{\delta }\int_{0}^{\infty }\frac{1}{x^{2+\delta }}f^{2}(x)dx\; .
\label{ineq}
\end{eqnarray}
\end{lemma}

\Proof  First, we recall the following Parseval identity
for Mellin transforms:
\begin{eqnarray*}
-\int_{0}^{\infty }\frac{f_{x}(x)(Hf)(x)}{x^{1+\delta }}dx=-\frac{1}{2\pi }%
\int_{-\infty }^{\infty }\overline{A(\lambda )}B(\lambda )d\lambda
\equiv I\;,
\end{eqnarray*}
with
\begin{eqnarray*}
A(\lambda )&=&\int_{0}^{\infty }x^{i\lambda -\frac{1}{2}-\frac{\delta }{2}%
}f_{x}(x)dx\;,
\\
B(\lambda )&=&\int_{0}^{\infty }x^{i\lambda -\frac{3}{2}-\frac{\delta }{2}%
}(Hf)(x)dx\;.
\end{eqnarray*}
Integration by parts in $A(\lambda )$ yields
\begin{eqnarray*}
A(\lambda )=-(i\lambda -\frac{1}{2}-\frac{\delta
}{2})\int_{0}^{\infty }x^{i\lambda -\frac{3}{2}-\frac{\delta
}{2}}f(x)dx\;.
\end{eqnarray*}
With respect to $B(\lambda )$ we can write
\begin{eqnarray*}
B(\lambda )&=&\int_{0}^{\infty }x^{i\lambda -\frac{3}{2}-\frac{\delta }{2}}%
\left[ \frac{1}{\pi }P.V.\int_{-\infty }^{+\infty }\frac{f(\xi )}{x-\xi }%
d\xi \right] dx
\\
&=&\int_{0}^{\infty }x^{i\lambda -\frac{3}{2}-\frac{\delta }{2}}\left[ \frac{1%
}{\pi }P.V.\int_{-\infty }^{0}\frac{f(\xi )}{x-\xi }d\xi +\frac{1}{\pi }%
P.V.\int_{0}^{\infty }\frac{f(\xi )}{x-\xi }d\xi \right] dx
\\
&=&\int_{0}^{\infty }x^{i\lambda -\frac{3}{2}-\frac{\delta }{2}}\left[ \frac{1%
}{\pi }\int_{0}^{\infty }\frac{f(\xi )}{x+\xi }d\xi +\frac{1}{\pi }%
P.V.\int_{0}^{\infty }\frac{f(\xi )}{x-\xi }d\xi \right] dx
\\
&=&\int_{0}^{\infty }x^{i\lambda -\frac{3}{2}-\frac{\delta }{2}}\left[ \frac{-x%
}{\pi }\int_{0}^{\infty }\frac{f(\xi )/\xi}{x+\xi }d\xi +\frac{x}{\pi }%
P.V.\int_{0}^{\infty }\frac{f(\xi )/\xi}{x-\xi }d\xi \right] dx
\\
&=&\int_{0}^{\infty }\left[- \frac{1}{\pi }\int_{0}^{\infty }\frac{x^{i\lambda -%
\frac{1}{2}-\frac{\delta }{2}}}{x+\xi }dx+\frac{1}{\pi
}P.V.\int_{0}^{\infty }\frac{x^{i\lambda -\frac{1}{2}-\frac{\delta
}{2}}}{x-\xi }dx\right] \frac{f(\xi)}{\xi}d\xi
\end{eqnarray*}
where we have used Fubini's theorem in order to exchange the order
of integration in $x$ and $\xi $. Using elementary complex
variable theory one can write
\begin{eqnarray*}
&&-\frac{1}{\pi }\int_{0}^{\infty }\frac{x^{i\lambda
-\frac{1}{2}-\frac{\delta
}{2}}}{x+\xi }dx+\frac{1}{\pi }P.V.\int_{0}^{\infty }\frac{x^{i\lambda -%
\frac{1}{2}-\frac{\delta }{2}}}{x-\xi }dx
\\
&=&\lim_{\substack{ R\rightarrow \infty  \\ \varepsilon \rightarrow
0}}\left[-\frac{1}{\pi }\frac{1}{1-e^{2\pi i(i\lambda -\frac{1}{2}-\frac{\delta }{2})}}%
\int_{\Gamma _{1}}\frac{z^{i\lambda -\frac{1}{2}-\frac{\delta }{2}}}{z+\xi }%
dz\right.
\\
&&\left. +\frac{1}{\pi }\frac{1}{1-e^{2\pi i(i\lambda -\frac{1}{2}-\frac{%
\delta }{2})}}\int_{\Gamma _{2}\backslash \left\{ c_{1},c_{2}\right\} }\frac{%
z^{i\lambda -\frac{1}{2}-\frac{\delta }{2}}}{z-\xi }dz\right]
\\
&\equiv& I_{1}+I_{2}
\end{eqnarray*}
where $\Gamma _{1}$ and $\Gamma _{2}$ are the paths in the complex
plane represented in Figures~1 and 2 respectively. Standard pole
integration for $I_{1}$ and the fact that
$\int_{\Gamma_{2}\backslash \left\{c_{1},c_{2}\right\}
}=-\int_{\left\{ c_{1},c_{2}\right\} }$ in $I_{2}$ (cf.\ Lemmas 2.2
and 2.3 in \cite{FV} where these integrals had to be computed for
a completely different purpose, for instance) yield then
\begin{equation*}
I_{1}+I_{2}=\left[ -\frac{1}{\sin \left( (-i\lambda
+\frac{1}{2}+\frac{\delta }{2})\pi \right) }+\cot \left( (-i\lambda
+\frac{1}{2}+\frac{\delta }{2})\pi \right) \right] \xi ^{i\lambda
-\frac{1}{2}-\frac{\delta }{2}}\;.
\end{equation*}
\begin{figure}[bp]
\begin{center}
\epsfig{file=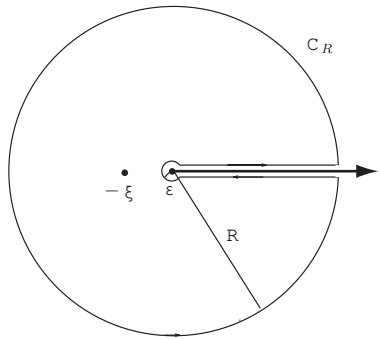} 
\end{center}
\vglue-18pt
\caption{Integration contour $\Gamma _{1}$.} \label{figure1}
\end{figure}
\begin{figure}[htbp]
\vglue-.75in
\begin{center}
\epsfig{file=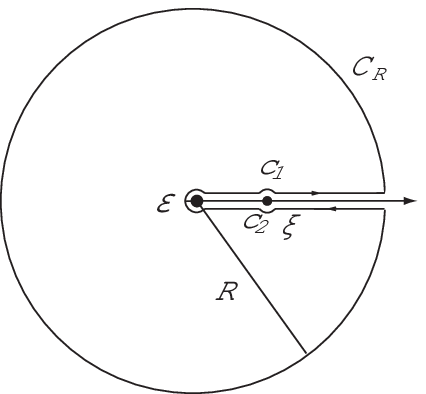} 
\end{center}
\vglue-18pt
\caption{Integration contour $\Gamma _{2}$.} \label{figure2}
\end{figure}

Hence
\begin{equation*}
B(\lambda )=\frac{-1+\cos \left( (-i\lambda +\frac{1}{2}+\frac{%
\delta }{2})\pi \right) }{\sin \left( (-i\lambda +\frac{1}{2}+\frac{\delta }{%
2})\pi \right) }F(\lambda )
\end{equation*}
with
\begin{equation*}
F(\lambda )\equiv \int_{0}^{\infty }\xi ^{i\lambda
-\frac{3}{2}-\frac{\delta }{2}}f(\xi )d\xi \; \pagebreak
\end{equation*}
and
\begin{eqnarray*}
I&=&\frac{1}{2\pi }\int_{-\infty }^{\infty }\frac{1-\cos \left( (-i\lambda +%
\frac{1}{2}+\frac{\delta }{2})\pi \right) }{\sin \left( (-i\lambda +\frac{1}{%
2}+\frac{\delta }{2})\pi \right) }(i\lambda +\frac{1}{2}+\frac{\delta }{2}%
)\left| F(\lambda )\right| ^{2}d\lambda
\\
&\equiv &\frac{1}{2\pi ^{2}}\int_{-\infty }^{\infty }M(\lambda
)\left|F(\lambda )\right| ^{2}d\lambda \;.
\end{eqnarray*}

In order to analyze $M(\lambda )$ we define now
\begin{eqnarray*}
z\equiv a+bi\;,\;a\equiv \left( \frac{1}{2}+\frac{\delta
}{2}\right) \pi \;,\;b\equiv \lambda \pi
\end{eqnarray*}
which implies, after some straightforward but lengthy
computations,
\begin{eqnarray}
M(\lambda )=z\frac{1-\cos \overline{z}}{\sin
\overline{z}}=\frac{a\sin a+b\sinh b}{\cosh b+\cos
a}+\frac{-a\sinh b+b\sin a}{\cosh b+\cos a}i \; . \label{ml}
\end{eqnarray}
Since $\left| F(\lambda )\right| ^{2}$ is symmetric in $\lambda $
and the imaginary part of $M(\lambda )$ is antisymmetric, 
\begin{eqnarray*}
I=\frac{1}{2\pi ^{2}}\int_{-\infty }^{\infty }{\rm Re}\left\{ M(\lambda
)\right\} \left| F(\lambda )\right| ^{2}d\lambda \;.
\end{eqnarray*}
Notice now from (\ref{ml}) that
\begin{eqnarray*}
\frac{1}{C}(1+\left| \lambda \right| )\leq {\rm Re}\left\{ M(\lambda
)\right\} \leq C(1+\left| \lambda \right| )
\end{eqnarray*}
so that
\begin{eqnarray*}
I\geq \frac{1}{2\pi C}\int_{-\infty }^{\infty }\left| F(\lambda
)\right|
^{2}d\lambda \geq C_{\delta }\int_{0}^{\infty }\frac{1}{x^{2+\delta }}%
f^{2}(x)dx
\end{eqnarray*}
where we have used the Plancherel  identity for Mellin
transforms:
\begin{eqnarray*}
\int_{0}^{\infty }\frac{1}{x^{2+\delta }}f^{2}(x)dx=\frac{1}{2\pi }%
\int_{-\infty }^{\infty }\left| F(\lambda )\right| ^{2}d\lambda\; .
\end{eqnarray*}
This completes the proof of the lemma.
\hfill\qed

\begin{rem}
Inequality (\ref{ineq}) can be extended by density to the restriction to $%
\mathbb{R}^{+}$ of any symmetric $f\in C^{1+\varepsilon
}(\mathbb{R})$ vanishing at the origin.
\end{rem}

In order to complete our blow-up argument, we have, from Cauchy's
inequality,
\begin{eqnarray*}
\int_{0}^{L}\frac{(1-\theta )}{x^{1+\delta }}dx&\leq& \left( \int_{0}^{L}\frac{%
(1-\theta )^{2}}{x^{2+\delta }}dx\right) ^{\frac{1}{2}}\left( \int_{0}^{L}%
\frac{1}{x^{\delta }}dx\right) ^{\frac{1}{2}}
\\
&\leq& \left( \frac{L^{1-\delta }}{1-\delta }\right)
^{\frac{1}{2}}\left(
\int_{0}^{\infty }\frac{(1-\theta )^{2}}{x^{2+\delta }}dx\right) ^{\frac{1}{2%
}}
\end{eqnarray*}
so that
\begin{eqnarray}
\int_{0}^{\infty }\frac{(1-\theta )^{2}}{x^{2+\delta }}dx\geq
C_{L,\delta }\left( \int_{0}^{L}\frac{(1-\theta )}{x^{1+\delta
}}dx\right) ^{2}\;. \label{ff}
\end{eqnarray}
From (\ref{eq3}), (\ref{ineq}) and (\ref{ff}) we deduce
\begin{eqnarray*}
\frac{d}{dt}\int_{0}^{L}\frac{(1-\theta )}{x^{1+\delta }}dx\geq
C_{L,\delta }\left( \int_{0}^{L}\frac{(1-\theta )}{x^{1+\delta
}}dx\right) ^{2}
\end{eqnarray*}
which yields a blow-up for
\begin{eqnarray*}
J\equiv \int_{0}^{L}\frac{(1-\theta )}{x^{1+\delta }}dx
\end{eqnarray*}
at finite time. Since
\begin{eqnarray*}
J\leq \int_{0}^{L}\frac{(1-\theta )}{x^{1+\delta }}dx\leq \sup_{x}\frac{%
1-\theta }{x}\int_{0}^{L}\frac{dx}{x^{\delta }}\leq \frac{L^{1-\delta }}{%
1-\delta }\sup_{x}\left| \theta _{x}\right|
\end{eqnarray*}
we conclude that $\left\| \theta _{x}\right\| _{L^{\infty }}$ must
blow up at finite time. This completes the proof of Theorem 2.1. \hfill\qed

\begin{rem}
In fact, numerical simulation by Morlet (see \cite{M1}) and
additional numerical experiments performed by ourselves (see
Figures 3 and 4) indicate that blow-up occurs at the maximum of
$\theta $ and   is such that a cusp develops at this point in
finite time.
\end{rem}

The figures below represent the profiles $\theta_x(x,t)$ and
$\theta(x,t)$ with initial data
\begin{eqnarray*}
\theta_0(x)=\begin{cases}
            (1 - x^2)^2, &{\rm if}  -1\leq x\leq 1 \\
            0,           &{\rm otherwise}
            \end{cases}
\end{eqnarray*}
at nine consecutive times.
\begin{figure}[htp]
 \begin{minipage}[t]{0.45\textwidth}
   \centering
\epsfig{file=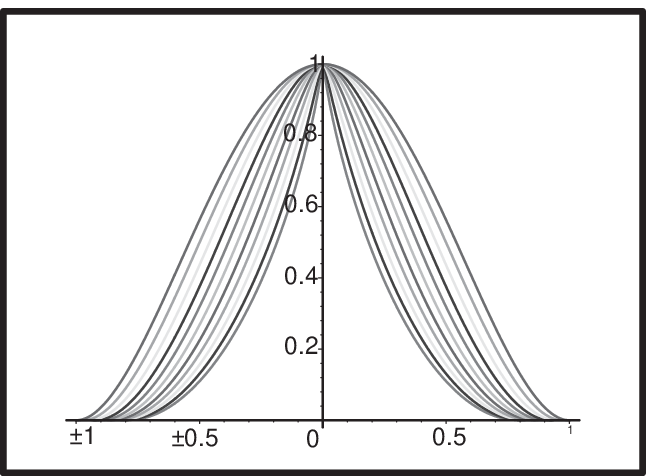}  
   \caption{$\theta(x,t)$}  \label{figure:3}
 \end{minipage}
\hspace{0.05\textwidth}%
 \begin{minipage}[t]{0.45\textwidth}
   \centering 
\epsfig{file=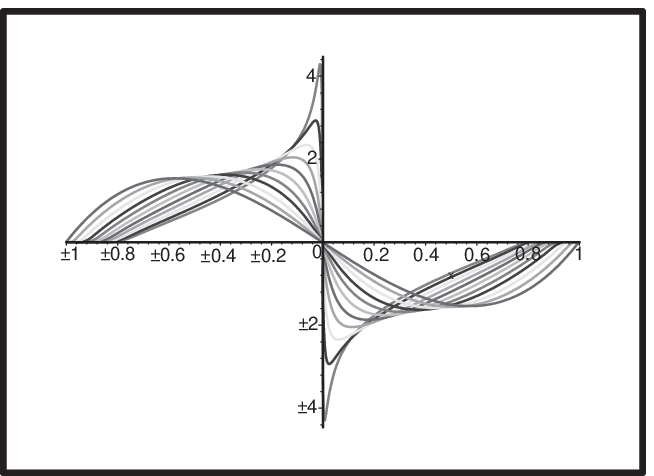} 
   \caption{$\theta_x(x,t)$}  \label{figure:4}
 \end{minipage}
\end{figure}

\section{The effect of viscosity}
Below we study the effect of viscosity ($\nu>0$) on the solutions
of (\ref{1eq}) with positive initial datum. First  

\begin{lemma}
Let $\theta$ be a $C^1$ solution of {\rm (\ref{1eq})} in $0\leq t \leq
T${\rm ,} with a nonnegative initial datum  $\theta_0 \in H^2(R)$. Then{\rm ,}
\begin{eqnarray}
 \qquad&1)& \qquad 0\leq \theta(x,t)\leq \|\theta_0\|_{L^{\infty}}\;  , \label{lm0}\\
&2)&\qquad  \|\theta\|_{L^{1}}(t)\leq
\|\theta_0\|_{L^{1}}\; , \label{lm1}\\
&3)&\qquad   \|\theta\|_{L^{2}}(t)\leq
\|\theta_0\|_{L^{2}}\quad\quad\text{and}\quad\quad
\int_{0}^{T}\|\Lambda^{\frac{\alpha}{2}}\theta\|^2_{L^{2}}dt\leq
\frac{\|\theta_0\|^2_{L^{2}}}{2\nu}\; . \label{lm2}
\end{eqnarray}
\end{lemma}

\Proof Since
\begin{eqnarray*}
\frac{d}{dt}\int|\Delta\theta|^2 dx = \int \Delta\theta \Delta(
H(\theta)\theta_x) dx\leq C\|\Delta\theta\|^3_{L^{2}}
\end{eqnarray*}
we have local solvability up to a time
$T=T(\|\theta_0\|_{H^{2}(R)})>0$ (without any restriction upon
the sign of $\theta_0$). Let us also observe that the same result
is true for the periodic version of (\ref{1eq}): $-\pi\leq x \leq
\pi$,
\begin{eqnarray*}
Hf(x)=P.V.\frac{1}{2\pi}\int_{-\pi }^{\pi
}\frac{f(x-y)}{{\rm tan}\frac{y}{2}}dy.
\end{eqnarray*}

We shall prove (\ref{lm0}) first in the periodic case: Let us
define $M(t)\equiv {\rm max}_x \theta(x,t)$, $m(t)\equiv
{\rm min}_ x\theta(x,t)$. It follows from the  H. Rademacher theorem that the
continuous Lipschitz functions $M(t)$, $m(t)$, admit ordinary
derivatives at almost every point t. Then we may argue as in
references \cite{CC} and \cite{CCF} to conclude that, at each
point of differentiability,  $M'(t)\leq 0$ and $m'(t)\geq
0$, implying (\ref{lm0}).

Let $\phi\in C^{\infty}_{0}(R)$ be such that $\phi\geq 0$,
$\phi(x)\equiv 1$ in $|x|\leq 1$ and $\phi(x)\equiv 0$ when
$|x|\geq 2$. With $R>0$ let us consider
$\theta_0^R(x)=\phi(\frac{x}{R})\theta_0(x)$ and let
$\theta^R(x,t)$ be the solution of the periodic problem
(\ref{1eq}) with initial data $\theta_0^R$ in $-\pi R\leq x \leq
\pi R$, $0\leq t \leq T = T(\theta_0)$.

We have that $0\leq \theta^R(x,t)\leq \|\theta_0\|_{L^{\infty}}$
with uniform estimates for $\nabla_x\theta^R$,
$\frac{\partial}{\partial t}\theta^R$. By compactness, we obtain a
sequence $\theta^{R_j}$, $R_j\rightarrow\infty$, converging
uniformly on compact sets to $\theta$, the solution of (\ref{1eq})
with initial data $\theta_0$. Then estimate (\ref{lm0}) follows.

To obtain  inequality (\ref{lm1}) we  proceed as follows:
\begin{eqnarray*}
\frac{d}{dt}\int\theta dx = \int H\theta\theta_x dx =
-\int\theta\Lambda\theta dx =
-\|\Lambda^{\frac12}\theta\|^2_{L^2},
\end{eqnarray*}
because $\int\Lambda^{\alpha}\theta dx= 0$.

 Next, observe that
\begin{eqnarray*}
\frac12\frac{d}{dt}\int\theta^2dx &=& \int\theta\theta_{x}H\theta
dx - \nu\int\theta\Lambda^{\alpha}\theta dx \\
&=& - \frac12\int\theta^{2}\Lambda\theta dx -
\nu\int|\Lambda^{\frac{\alpha}{2}}\theta|^2dx.
\end{eqnarray*}
On the other hand
\begin{eqnarray*}
\int\theta^{2}\Lambda\theta dx = \int\int\frac{[\theta(x) +
\theta(y)]}{2}\frac{(\theta(x) - \theta(y))^2}{(x - y)^2}dxdy\geq
0
\end{eqnarray*}
and the proof of the third part of the lemma follows.
\Subsec{Global existence with $\alpha>1$}

\begin{thm}
Let $0\leq\theta_0\in H^2(R)${\rm ,} $\nu>0$ and $\alpha>1$. Then there
exists a constant $C${\rm ,} depending only on $\theta_0$ and $\nu${\rm ,}
such that for $t\geq 0$\/{\rm :}\/
\begin{eqnarray}
&1)&\qquad \|\Lambda^{\frac{1}{2}}\theta\|_{L^{2}}(t)\leq C\; , \label{lm3}\\
&2)&\qquad \|\Lambda\theta\|_{L^{2}}(t)\leq C(1 + t)\; ,\label{lm4}\\
&3)&\qquad \|\Delta\theta\|_{L^{2}}(t)\leq Ce^{Ct^3}\; .\label{lm5}
\end{eqnarray}
\end{thm}

\Proof  Integration by parts and the formula for the
Hilbert transform
\begin{eqnarray*}
2H(fH(f))= (H(f))^2 - f^2
\end{eqnarray*}
 yield 
\begin{eqnarray} &&\label{1/2}\\
\frac12\frac{d}{dt}\int|\Lambda^{\frac{1}{2}}\theta|^2dx &=&
\int\Lambda\theta\theta_{x}H\theta
dx  - \nu\int|\Lambda^{\frac 12 + \frac{\alpha}{2}}\theta|^2dx\nonumber\\
&=& - \int\theta H(\theta_{x}H\theta_{x})
dx  - \nu\int|\Lambda^{\frac12 + \frac{\alpha}{2}}\theta|^2dx\nonumber\\
&=& - \frac12\int\theta (H\theta_{x})^2dx
  + \frac12\int\theta (\theta_{x})^2dx
  - \nu\int|\Lambda^{\frac12 + \frac{\alpha}{2}}\theta|^2dx\nonumber\\
  &\leq& \|\theta_0\|_{L^{\infty}}\|\Lambda\theta\|_{L^{2}}^2 -
\nu\int|\Lambda^{\frac12 + \frac{\alpha}{2}}\theta|^2dx.
\nonumber
\end{eqnarray}
Since
$$
\|\Lambda\theta\|_{L^{2}}^2\leq R^{2-\alpha}\|\Lambda^{
\frac{\alpha}{2}}\theta\|_{L^{2}}^2 + R^{1 -
\alpha}\|\Lambda^{\frac12 + \frac{\alpha}{2}}\theta\|_{L^{2}}^2,
$$
  by taking $R$ sufficiently large and applying inequality
(\ref{lm2}), we obtain the desired inequality (\ref{lm3}).

Applying $\Lambda$ operator to both sides of equation (\ref{1eq}),
multiplying by $\Lambda\theta$ and integrating in $x$, we obtain
\begin{eqnarray} \qquad 
\frac12\frac{d}{dt}\int|\Lambda\theta|^2dx &=&
\int\Lambda\theta\Lambda(\theta_{x}H\theta) dx  - \nu\|\Lambda^{1
+ \frac{\alpha}{2}}\theta\|_{L^{2}}^2\label{1}\\
&=& -\frac12\int(\theta_x)^2\Lambda\theta dx - \nu\|\Lambda^{1 +
\frac{\alpha}{2}}\theta\|_{L^{2}}^2\nonumber\\
&\leq& C\int|\Lambda\theta|^3dx- \nu\|\Lambda^{1 +
\frac{\alpha}{2}}\theta\|_{L^{2}}^2\; , \nonumber
\end{eqnarray}
where we have used the isometry of Hilbert transform in $L^2$,
 integration by parts and finally Cauchy's inequality together with the boundedness of
 Hilbert transform in $L^3$.

In order to estimate $\|\Lambda\theta\|_{L^{3}}$ we make use of
Hausdorff-Young's inequality
\begin{eqnarray}
\|\Lambda\theta\|_{L^{3}}\leq
\|\widehat{\Lambda\theta}\|_{L^{\frac32}} =
\left(\int|\xi|^{\frac32}|\widehat{\theta}(\xi)|^{\frac32}d\xi\right)^{\frac23}.\label{L3}
\end{eqnarray}

Picking now $\bar{\alpha}\in (1, \alpha)$ and using Cauchy's
inequality we obtain
\begin{eqnarray*}
\left(\int|\xi|^{\frac32}|\widehat{\theta}(\xi)|^{\frac32}d\xi\right)^{\frac23}
&\leq& \left(\int|\xi|^{2 + \bar{\alpha}
}|\widehat{\theta}(\xi)|^{2}d\xi\right)^{\frac13}
\left(\int|\xi|^{1 - \bar{\alpha}
}|\widehat{\theta}(\xi)|d\xi\right)^{\frac13}\\
&\equiv& I_1^{\frac13}\cdot I_2^{\frac13}.
\end{eqnarray*}

For $I_1$ we get
\begin{eqnarray}
I_1&=& \int_{|\xi|\leq R}|\xi|^{2 + \bar{\alpha}
}|\widehat{\theta}(\xi)|^{2}d\xi + \int_{|\xi|\geq R}|\xi|^{2 +
\bar{\alpha} }|\widehat{\theta}(\xi)|^{2}d\xi\label{I1}\\
&\leq& R^{2 + \bar{\alpha}}\|\theta\|_{L^{2}}^2 +
\frac{1}{R^{\alpha-\bar{\alpha}}}\int_{|\xi|\geq R}|\xi|^{2 +
\alpha }|\widehat{\theta}(\xi)|^{2}d\xi\nonumber\\
&\leq& R^{2 + \bar{\alpha}}\|\theta\|_{L^{2}}^2 +
\frac{1}{R^{\alpha-\bar{\alpha}}}\|\Lambda^{1 +
\frac{\alpha}{2}}\theta\|_{L^{2}}^2.\nonumber
\end{eqnarray}

With respect to $I_2$ one can estimate
\begin{eqnarray} \qquad 
I_2&=& \int_{|\xi|\leq 1}|\xi|^{1 - \bar{\alpha}
}|\widehat{\theta}(\xi)|d\xi + \int_{|\xi|\geq 1}|\xi|^{1 -
\bar{\alpha} }|\widehat{\theta}(\xi)|d\xi\label{I2}\\
& \leq& \int_{|\xi|\leq 1}|\widehat{\theta}(\xi)|d\xi +
\int_{|\xi|\geq 1}|\xi|^{\frac12
}|\widehat{\theta}(\xi)||\xi|^{\frac12 - \bar{\alpha} }d\xi\nonumber\\
& \leq& \int_{|\xi|\leq 1}\|\theta\|_{L^{1}}d\xi +
\left(\int_{|\xi|\geq 1}|\xi|^{1- 2\bar{\alpha}
}d\xi\right)^{\frac12}\left(\int_{|\xi|\geq
1}|\xi||\widehat{\theta}(\xi)|^2d\xi\right)^{\frac12}\nonumber\\
&\leq& \|\theta\|_{L^{1}} +
c_{\alpha}\|\Lambda^{\frac12}\theta\|_{L^{2}}\; .\nonumber
\end{eqnarray}

From (\ref{I2}), (\ref{lm1}) and (\ref{lm3}) it follows that
\begin{eqnarray}
I_2\leq C.\label{I22}
\end{eqnarray}
Hence by (\ref{L3}), (\ref{I1}) and (\ref{I22}) we get
\begin{eqnarray}
\|\Lambda\theta\|_{L^{3}}\leq C^{\frac13}(R^{\frac{2 +
\bar{\alpha}}{3}}\|\theta\|_{L^{2}}^{\frac23} +
\frac{1}{R^{\frac{\alpha-\bar{\alpha}}{3}}}\|\Lambda^{1 +
\frac{\alpha}{2}}\theta\|_{L^{2}}^{\frac23})\; . \label{lll}
\end{eqnarray}
To finish let us take $R$ sufficiently large together with
(\ref{1}), (\ref{lll}) and (\ref{lm2}) to conclude that
\begin{eqnarray*}
\frac12\frac{d}{dt}\int|\Lambda\theta|^2dx \leq
C\|\theta_0\|_{L^{2}}^2
\end{eqnarray*}
from which  (\ref{lm4})  follows.

Finally let us consider
\begin{eqnarray}\quad 
\frac12\frac{d}{dt}\|\Delta\theta\|_{L^2}^2&=&
\int\Delta\theta\Delta(\theta_xH\theta)dx -
\nu\|\Lambda^{2+\frac{\alpha}{2}}\theta\|_{L^2}^2\label{l1}\\
&\leq& C
\|\Delta\theta\|_{L^{\infty}}\|\Delta\theta\|_{L^{2}}\|\Lambda\theta\|_{L^{2}}
- \nu\|\Lambda^{2+\frac{\alpha}{2}}\theta\|_{L^2}^2\nonumber
\end{eqnarray}
and let us observe that
\begin{eqnarray}
\|\Delta\theta\|^2_{L^{\infty}}\leq C (
\|\Lambda^{2+\frac{\alpha}{2}}\theta\|^2_{L^2} +
\|\theta\|^2_{L^2}).\label{l2}
\end{eqnarray}
Therefore, by Holder's inequality,
\begin{eqnarray*}
\|\Delta\theta\|_{L^{\infty}}\|\Delta\theta\|_{L^{2}}\|\Lambda\theta\|_{L^{2}}\leq
\frac{\delta}{2}\|\Delta\theta\|^2_{L^{\infty}} +
\frac{1}{2\delta}\|\Delta\theta\|^2_{L^{2}}\|\Lambda\theta\|_{L^{2}}^2,
\end{eqnarray*}
and inequality (\ref{l2}) we estimate the first term at the 
right-hand side of (\ref{l1}), and conclude  that choosing $\delta$ small
enough,
\begin{eqnarray*}
\frac{d}{dt}\|\Delta\theta\|_{L^2}^2\leq
C(\|\Lambda\theta\|^2_{L^{2}}\|\Delta\theta\|_{L^2}^2 +
\|\theta\|^2_{L^2})
\end{eqnarray*}
which implies the estimate
\begin{eqnarray*}
\|\Delta\theta\|_{L^2}^2&\leq& \|\Delta\theta_0\|_{L^2}^2
e^{C\int^{t}_{0}\|\Lambda\theta\|^2_{L^{2}}ds} +
C\int_0^t\|\theta\|^2_{L^2}
e^{\int_s^t\|\Lambda\theta\|^2_{L^{2}}d\sigma}ds.
\end{eqnarray*}
By (\ref{lm2}) and (\ref{lm4}),  (\ref{lm5})  then follows  for
some large enough $C$.

\Subsec{Small data results for $\alpha=1$}
In the critical case $\alpha = 1$ we have the following global
existence result for small data.
\begin{thm}
Let $\nu >0${\rm ,} $\alpha = 1${\rm ,} $0\leq\theta_0\in H^1$ and assume that
the initial data satisfy $\|\theta_0\|_{L^{\infty}}< \nu $. Then
there exists a unique solution to {\rm (\ref{1eq})} which belongs to
$H^1$ for all time $t>0$.
\end{thm}

\Proof  From the previous inequality (\ref{1/2}) we have for
$\alpha = 1$
\begin{eqnarray}
\frac12\frac{d}{dt}\int|\Lambda^{\frac{1}{2}}\theta|^2dx &\leq&
\|\theta_0\|_{L^{\infty}}\|\Lambda\theta\|_{L^{2}}^2 -
\nu\int|\Lambda\theta|^2dx\label{312eq}\\
&=& (\|\theta_0\|_{L^{\infty}} - \nu)\|\Lambda\theta\|_{L^{2}}^2,
\nonumber
\end{eqnarray}
which implies that if $\|\theta_0\|_{L^{\infty}}< \nu $, then
\begin{eqnarray}
\|\Lambda^\frac12\theta\|_{L^2}(t)\leq
\|\Lambda^\frac12\theta_0\|_{L^2}\quad and\quad
\int_0^t\|\Lambda\theta\|_{L^2}^2ds\leq
C\|\Lambda^\frac12\theta_0\|^2_{L^2}.\label{l3}
\end{eqnarray}

From (\ref{1}) we get
\begin{eqnarray}
\frac12\frac{d}{dt}\int|\Lambda\theta|^2dx &\leq&
\frac12\int|\Lambda\theta|^3dx-
\nu\|\Lambda^{\frac32}\theta\|_{L^{2}}^2.\label{312eq}
\end{eqnarray}

Since
\begin{eqnarray*}
\|\Lambda\theta\|_{L^{3}}^3 \leq \|\Lambda\theta\|_{L^{2}}^2\cdot
\|\Lambda\theta\|_{{\rm BMO}}
\end{eqnarray*}
and
\begin{eqnarray*}
\|\Lambda\theta\|_{{\rm BMO}}&\leq& C\|\Lambda^{\frac32}\theta\|_{L^{2}}
 \end{eqnarray*}
(we refer to \cite{stein} for the corresponding definitions and
properties of the functions of bounded mean oscillation ({\rm BMO})), we
obtain
\begin{eqnarray*}
\frac12\frac{d}{dt}\int|\Lambda\theta|^2dx &\leq&
C\|\Lambda\theta\|_{L^{2}}^2\|\Lambda^{\frac32}\theta\|_{L^{2}} -
\nu\|\Lambda^{\frac32}\theta\|_{L^{2}}^2\\
&\leq&\frac{C^2}{4\nu}\|\Lambda\theta\|_{L^{2}}^4.
\end{eqnarray*}
Together with inequalities (\ref{l3}) this  allows us to complete the
proof of the theorem.


\references {999}
\bibitem[1]{B} \name{G. R. Baker, X. Li}, and \name{A. C. Morlet}, Analytic structure
of 1D-transport equations with nonlocal fluxes, {\it Physica D}
{\bf 91}  (1996), 349--375.

\bibitem[2]{BM} \name{A. L. Bertozzi} and \name{A. J. Majda}, {\it  
Vorticity and the Mathematical
Theory of Incompresible Fluid Flow\/}, Cambridge University Press,
Cambridge, 2002.



\bibitem[3]{CLM}  \name{P. Constantin, P. Lax}, and \name{A. Majda}, A simple one-dimensional
model for the three dimensional vorticity, \textit{Comm.\ Pure
Appl.\ Math.\/} \textbf{38} (1985), 715--724.

\bibitem[4]{CC}  \name{A. C\'{o}rdoba} and \name{D. C\'{o}rdoba}, A maximum principle applied to
quasi-geostrophic equations, \textit{Comm.\ Math.\ Phys.\/} \textbf{249}
(2004), 511--528.

\bibitem[5]{CCF}  \name{D. Chae, A. C\'ordoba, D. 
C\'ordoba}, and \name{M. A. Fontelos}, Finite time
singularities in a 1D model of the quasi-geostrophic equation.
\textit{Adv.\ Math.\/} \textbf{194} (2005), 203--223.

\bibitem[6]{De1} \name{S. De Gregorio}, On a one-dimensional model for the
three-dimensional vorticity equation, \textit{J. Statist.\ Phys.\/}
{\bf 59} (1990), 1251--1263.

\bibitem[7]{De2} \bibline, A partial differential equation
arising in a 1D model for the 3D vorticity equation, \textit{Math.\
Methods Appl.\ Sci.\/} {\bf 19} (1996), 1233--1255.


\bibitem[8]{FV}  \name{M. A. Fontelos} and \name{J. J. L. Vel\'{a}zquez}, A free boundary
problem for the Stokes system with contact lines, \textit{Comm.\
Partial Differential Equations} {\bf 23} (1998), 
1209--1303.



\bibitem[9]{Moore} \name{D. W. Moore}, The spontaneous appearance of a
singularity in the shape of an evolving vortex sheet,
\textit{Proc.\ R. Soc.\ London} A {\bf 365} (1979),
105-119.

\bibitem[10]{M} \name{A. Morlet}, Further properties of a continuum of model
equations with globally defined flux. \textit{J. Math. Anal.
Appl.} \textbf{221} (1998), 132-160.

\bibitem[11]{M1} \bibline, Some further results for a one-dimensional
transport equation with nonlocal flux, \textit{Comm.\ Appl.\ Anal.\/}
{\bf 1} (1997), 315--336.


\bibitem[12]{Sakajo} \name{T. Sakajo}, On global solutions for the
Constantin-Lax-Majda equation with a generalized viscosity term,
\textit{Nonlinearity} \textbf{16} (2003), 1319--1328.

\bibitem[13]{Schochet} \name{S. Schochet}, Explicit solutions of the viscous
model vorticity equation, \textit{Comm.\ Pure Appl.\ Math.\/}
\textbf{41} (1986), 531--537.



\bibitem[14]{stein} \name{E. Stein}, {\it Harmonic Analysis\/}: {\it Real Variable 
Methods\/}, {\it Orthogonality
and Oscillatory Integrals}, {\it Princeton Math.\ Series\/} {\bf 43}, Princeton Univ. Press, Princeton, NJ,
1993.

\bibitem[15]{Yang} \name{Y. Yang}, Behavior of solutions of model equations for incompressible fluid flow,
\textit{J. Differential Equations} \textbf{125} (1996), 133--153.

\bibitem[16]{Vasudeva} \name{M. Vasudeva}, The Constantin-Lax-Majda model vorticity equation revisited,
\textit{J. Indian Inst.\ Sci.\/} {\bf 78} (1998), 109--117.

\bibitem[17]{Vaswegert} \name{M. Vasudeva} and \name{E. Wegert}, Blow-up in a modified
Constantin-Lax-Majda model for the vorticity equation, \textit{Z.
Anal.\ Anwend.\/ } {\bf 18} (1999), 183--191.

\Endrefs
\end{document}